\documentclass[11pt]{article}
\usepackage{amsmath,amssymb,amsthm}
\usepackage[margin=2.5cm]{geometry}
\usepackage{enumitem}
\usepackage{hyperref}
\usepackage{authblk}

\newtheorem{theorem}{Theorem}[section]
\newtheorem{lemma}[theorem]{Lemma}

\theoremstyle{remark}
\newtheorem{remark}[theorem]{Remark}

\allowdisplaybreaks

\title{A New Proof of Wei's Theorem}
\author{Zhiqiang Wei \quad Mengyan Xing \quad Xingya Zhao}

\date{\today}

\begin{document}
\maketitle

\begin{abstract}
In the paper \emph{Structure and realizability for rational maps}, Wei
proved the following. Let \(k\) nontrivial partitions of an integer
\(d\ge 3\) be given, satisfying the Riemann--Hurwitz formula for the case of
maps between two spheres. If \(k>\ell+1\), where
\[
\ell=\min_i|\pi_i|
\]
is the shortest partition length, then there exists a rational map whose
branch datum is given by these partitions. In this paper we give a new proof of Wei's theorem on realizing a tuple of
partitions as the ramification data of a rational function. The proof is
inductive on the shortest partition length. Two auxiliary results are
established: a two-cut partition theorem and a gluing theorem for marked
spherical coverings. These reduce any counterexample with shortest length
\(\ell\ge3\) to auxiliary data with either fewer factors or smaller shortest
length, while preserving the total defect and the marked conditions. The base
case \(\ell\le2\) is quoted from the literature.
\end{abstract}
\maketitle
\section{Introduction}

Wei's theorem \cite[Theorem~1.2]{Wei} concerns the realization of a tuple of
partitions \(\pi_1,\ldots,\pi_k\vdash d\) as the ramification data of a
rational function. The purpose of this note is to give a new proof of this
theorem, based on two operations on partitions and spherical coverings:
\emph{two-cut partitions} and \emph{marked gluing}.

The proof is inductive on
\[
\ell=\min_i|\pi_i|.
\]
Assuming the base case \(\ell\le2\), we show that any counterexample with
\(\ell\ge3\) can be reduced to auxiliary data with either fewer factors or
smaller shortest length. The reduction uses a two-cut partition theorem
(Theorem~\ref{thmA}) and a marked gluing theorem for spherical coverings
(Theorem~\ref{thmB}), together with an additive-defect merging lemma
(Lemma~\ref{lem:2}). The final gluing uses common letters in the marked
cycles; Riemann--Hurwitz and the Riemann existence theorem then yield the
required rational function. Note that the realizability for the case $l=1$ is well known \cite{EKS84}. The realizability in the case $l=2$ is proved by Pakovich \cite{Pakovich}.

\section{Notation and the first lemma}

For a partition \(\pi\vdash n\), let \(|\pi|\) be its number of parts and
define the \emph{defect}
\[
\nu(\pi)=n-|\pi|.
\]
The partition \([1^n]\) is called \emph{trivial}. The symbol \(\sqcup\)
denotes multiset union. An \emph{allowed cut} selects several complete parts
of a partition and allows one further part to be shortened; if the
shortening is zero, no actual cut occurs.

\begin{lemma}[Integer interval of attainable defects]\label{lem:1}
Let \(\pi\vdash d\) and \(1\le b<d\). The defects of partitions
\(\mu\vdash b\) obtained by one allowed cut form an integer interval
\[
[m(\pi,b),M(\pi,b)]\cap\mathbb Z,
\]
and
\[
m(\pi,b)\le \left\lfloor\frac{b\nu(\pi)}{d}\right\rfloor,
\qquad
M(\pi,b)\ge \min\left\{\left\lfloor\frac b2\right\rfloor,\nu(\pi)\right\}.
\]
If \(\nu(\pi)>0\), then \(m(\pi,b)\le \nu(\pi)-1\).
\end{lemma}

\begin{proof}
Write \(r=|\pi|\). A partition \(\mu\) with \(|\mu|=q\) is attainable iff
the sum of the largest \(q\) parts of \(\pi\) is at least \(b\), and the
sum of the smallest \(q-1\) parts is at most \(b-1\). Necessity is clear.
For sufficiency, first take the smallest \(q\) parts. If their sum reaches
\(b\), shorten the largest of them. Otherwise, successively replace chosen
parts by larger ones until the total first reaches \(b\), then shorten the
newly chosen part. Before the last replacement the total is less than
\(b\), so the sum of the remaining \(q-1\) chosen parts is less than \(b\);
hence the shortened part is positive.

In the above criterion, the first condition is monotone increasing in
\(q\), while the second is monotone decreasing in \(q\). Hence the
attainable \(q\), and therefore the defects \(b-q\), form an integer
interval.

For the upper bound on \(m(\pi,b)\), take
\[
q=\left\lceil\frac{br}{d}\right\rceil.
\]
Then the sum of the largest \(q\) parts is at least
\[
\frac{qd}{r}\ge b,
\]
and the sum of the smallest \(q-1\) parts is at most
\[
\frac{(q-1)d}{r}<b.
\]
Thus this \(q\) is attainable, and consequently
\[
m(\pi,b)\le b-q
=b-\left\lceil\frac{br}{d}\right\rceil
=\left\lfloor\frac{b(d-r)}{d}\right\rfloor
=\left\lfloor\frac{b\nu(\pi)}{d}\right\rfloor,
\]
where we used \(r=d-\nu(\pi)\) and the identity
\(b-\lceil x\rceil=\lfloor b-x\rfloor\) for \(b\in\mathbb Z\).

For the lower bound on \(M(\pi,b)\), we distinguish two cases. If the sum
of all nontrivial parts of \(\pi\) is less than \(b\), take all of them
and fill up with \(1\)'s; the resulting defect is exactly \(\nu(\pi)\).
Otherwise, select nontrivial parts successively until the sum reaches
\(b\), and shorten the last selected part. Since each nontrivial part is
at least \(2\), at most
\[
\left\lceil\frac b2\right\rceil
\]
parts are needed. Hence the defect is at least
\[
b-\left\lceil\frac b2\right\rceil
=\left\lfloor\frac b2\right\rfloor.
\]
Therefore
\[
M(\pi,b)\ge \min\left\{\left\lfloor\frac b2\right\rfloor,\nu(\pi)\right\}.
\]

Finally, since \(b<d\) and \(\nu(\pi)>0\),
\[
\left\lfloor\frac{b\nu(\pi)}{d}\right\rfloor\le \nu(\pi)-1,
\]
so the first estimate gives
\[
m(\pi,b)\le \nu(\pi)-1.
\]
This completes the proof.
\end{proof}

\section{The two-cut partition theorem}

\begin{theorem}[Two cuts]\label{thmA}
Let \(\ell\ge3\) and \(\pi_1,\ldots,\pi_{\ell+1}\vdash d\) satisfy
\[
\nu_i:=\nu(\pi_i)\ge\ell,
\qquad
\sum_{i=1}^{\ell+1}\nu_i=d+\ell-2.
\]
Let \(b\ge1\) with \(\ell b\le d\) and put \(D=d-b\). Then there exist
\(\mu_i\vdash b\), \(\lambda_i\vdash D\), and distinct indices \(I,J\)
such that:
\begin{enumerate}[label=(\arabic*)]
\item if \(i\notin\{I,J\}\), then \(\pi_i=\mu_i\sqcup\lambda_i\);
\item if \(i\in\{I,J\}\), there exist marked parts \(u_i\in\mu_i\),
      \(v_i\in\lambda_i\) with
      \[
      \pi_i=(\mu_i\setminus\{u_i\})\sqcup
      (\lambda_i\setminus\{v_i\})\sqcup\{u_i+v_i\};
      \]
\item \(\sum_i\nu(\mu_i)=b-1\);
\item \(\nu(\lambda_i)>0\) for \(i\notin\{I,J\}\), and
      \[
      0<\nu(\lambda_I)+\nu(\lambda_J)\le D-2;
      \]
\item if \(b\ge2\) and \(u_I=u_J=1\), then
      \(\nu(\mu_I)+\nu(\mu_J)\le b-2\).
\end{enumerate}
\end{theorem}

\begin{proof}
Put \(S=d+\ell-2\). Any two distinct indices \(p,q\) satisfy
\begin{equation}\label{eq:2}
\nu_p+\nu_q\le S-(\ell-1)=d-1,
\end{equation}
since the remaining \(\ell-1\) rows have defect at least \(1\).

\emph{Case \(b=1\).} A partition without part \(1\) has defect at least
\(d/2\); by \eqref{eq:2} there is at most one such row. Choose it and one
further row as \(I,J\), split a part \(a\ge2\) in each into \(1+(a-1)\),
and take an existing \(1\) from every other row. Then all
\(\mu_i=[1]\), and
\[
0<\nu_I+\nu_J-2=\nu(\lambda_I)+\nu(\lambda_J)\le d-3=D-2.
\]
The remaining rows keep positive large-side defect.

\emph{Case \(b\ge2\).} Let \(e_i\) be the number of \(1\)'s in
\(\pi_i\). If \(e_i<b\), then \(\nu_i\ge(d-e_i)/2>(d-b)/2\). Since
\(d\ge3b\), three such rows would give
\[
S>\frac{3(d-b)}{2}+\ell-2\ge d+\ell-2,
\]
a contradiction. Hence at most two rows have \(e_i<b\), and for such rows
\begin{equation}\label{eq:3}
e_i<b\Longrightarrow \nu_i\ge b+1.
\end{equation}
Choose the two rows with \(e_i<b\) if there are two; otherwise choose one
such row together with a row of maximal defect; if there is none, choose
two rows of maximal defect. Call the chosen rows \(P,Q\). All other rows
contain at least \(b\) ones.

For \(i=P,Q\), set
\[
m_i=m(\pi_i,b),\qquad
M_i=M(\pi_i,b),\qquad
U_i=\min\{M_i,\nu_i-1\}.
\]
Lemma~\ref{lem:1} and \eqref{eq:2} give
\begin{equation}\label{eq:4}
m_P+m_Q
\le \left\lfloor\frac{b(\nu_P+\nu_Q)}{d}\right\rfloor
\le b-1.
\end{equation}
We claim
\begin{equation}\label{eq:5}
U_P+U_Q\ge b-1.
\end{equation}
If both \(\nu_P\) and \(\nu_Q\) are at least \(\lfloor b/2\rfloor+1\),
then Lemma~\ref{lem:1} gives \(U_i\ge\lfloor b/2\rfloor\) for
\(i=P,Q\), and \eqref{eq:5} follows. Otherwise, after swapping names,
let \(t=\nu_Q\le\lfloor b/2\rfloor\). We claim that \(\nu_i\le t\) for
all \(i\ne P\). Indeed, by \eqref{eq:3} and
\(t\le\lfloor b/2\rfloor<b+1\), the row \(Q\) cannot satisfy
\(e_Q<b\); hence \(e_Q\ge b\). If there is a row with \(e_i<b\), then
it must be \(P\), and \(Q\) was chosen as a row of maximal defect among
the remaining rows, so every \(i\ne P,Q\) satisfies \(\nu_i\le\nu_Q=t\).
If no row has \(e_i<b\), then \(P,Q\) are the two rows of maximal defect,
so again every \(i\ne P,Q\) satisfies \(\nu_i\le\nu_Q=t\). Thus
\(\nu_i\le t\) for all \(i\ne P\).

Then
\[
|\pi_P|=d-\nu_P=\sum_{i\ne P}\nu_i-\ell+2
\le \ell t-\ell+2<\ell t.
\]
If \(t\le|\pi_P|\), the largest \(t\) parts of \(\pi_P\) sum to at least
\(td/|\pi_P|>d/\ell\ge b\); if \(t>|\pi_P|\), all parts of \(\pi_P\)
already sum to at least \(b\). Hence at most \(t\) parts suffice to reach
\(b\), so \(M_P\ge b-t\). Also
\[
\nu_P\ge d+\ell-2-\ell t\ge \ell(b-t)+\ell-2\ge b-t+1,
\]
so \(U_P\ge b-t\). On the other hand, Lemma~\ref{lem:1} gives
\[
M_Q\ge\min\left\{\left\lfloor\frac b2\right\rfloor,\nu_Q\right\}=t,
\]
hence \(U_Q\ge t-1\). Thus \eqref{eq:5} holds.

By \eqref{eq:4}, \eqref{eq:5} and the interval property, choose
\[
s_P\in[m_P,U_P],\qquad
s_Q\in[m_Q,U_Q],\qquad
s_P+s_Q=b-1.
\]
Take the corresponding \(\mu_P,\mu_Q\), and set \(\mu_i=[1^b]\) for all
other rows. At this stage at most two actual cuts have been made and the
total small-side defect is \(b-1\).

It remains to normalize the cuts. If fewer than two cuts have been made,
choose the required number of rows outside \(P,Q\). Return one original
\(1\) from the small side to the large side, split some \(a\ge2\) on the
large side into \(1+(a-1)\), and put the new \(1\) on the small side. The
small partition remains \([1^b]\), while the large-side defect becomes
\(\nu_i-1\ge\ell-1>0\).

Finally, suppose two actual cuts have been made and both small-side marked parts are \(1\). Define
\[
c(\eta)=\sum_{a\ge 2}(a-2).
\]
Let \(e_i\) be the number of \(1\)'s in \(\pi_i\), and let \(r_i\) be the number of parts of \(\pi_i\) different from \(1\). Then \(|\pi_i|=e_i+r_i\), \(\nu_i=d-e_i-r_i\), and
\[
c(\pi_i)=\sum_{a\ge 2}(a-2)=\sum_{a\ge 2}a-2r_i=d-e_i-2r_i.
\]
Hence
\[
2\nu_i=2d-2e_i-2r_i=(d-e_i-2r_i)+(d-e_i)=c(\pi_i)+d-e_i,
\]
so
\[
2\nu_i=d-e_i+c(\pi_i). \tag{6}
\]
For any two distinct indices \(p,q\), (1) gives \(\nu_p+\nu_q\le d-1\). Taking \(p=P,q=Q\) and using (6), we obtain
\[
(d-e_P+c(\pi_P))+(d-e_Q+c(\pi_Q))=2(\nu_P+\nu_Q)\le 2d-2.
\]
After cancellation,
\[
e_P+e_Q\ge c(\pi_P)+c(\pi_Q)+2>c(\pi_P)+c(\pi_Q).
\]
This proves (5):
\[
e_P+e_Q>c(\pi_P)+c(\pi_Q). \tag{5}
\]

Now let \(f_i\) be the number of original \(1\)'s of \(\pi_i\) used on the small side of \(\mu_i\), excluding the marked \(1\)'s produced by cuts, and let \(r_i'\) be the number of non-\(1\) parts of \(\mu_i\). Since each cut row has exactly one marked \(1\) produced by the cut on the small side,
\[
|\mu_i|=f_i+1+r_i',\qquad \nu(\mu_i)=b-|\mu_i|=b-1-f_i-r_i'.
\]
On the other hand,
\[
c(\mu_i)=\sum_{a\ge 2}(a-2)=\sum_{a\ge 2}a-2r_i'=(b-f_i-1)-2r_i',
\]
so
\[
c(\mu_i)=\nu(\mu_i)-r_i'. \tag{7}
\]
Because the total small-side defect is \(\sum_i\nu(\mu_i)=b-1\), and the ordinary rows \(\mu_i=[1^b]\) have defect \(0\), we have
\[
\nu(\mu_P)+\nu(\mu_Q)=b-1.
\]
Combining this with (7),
\[
\begin{aligned}
f_P+f_Q
&=\bigl(b-1-\nu(\mu_P)-r_P'\bigr)
 +\bigl(b-1-\nu(\mu_Q)-r_Q'\bigr)\\
&=2b-2-\bigl(\nu(\mu_P)+\nu(\mu_Q)\bigr)-r_P'-r_Q'\\
&=b-1-r_P'-r_Q',
\end{aligned}
\]
while
\[
c(\mu_P)+c(\mu_Q)
=\nu(\mu_P)+\nu(\mu_Q)-r_P'-r_Q'
=b-1-r_P'-r_Q'.
\]
Thus
\[
f_P+f_Q=c(\mu_P)+c(\mu_Q). \tag{8}
\]
Since an allowed cut can only decrease or preserve the value of \(c\), we have
\[
c(\mu_P)+c(\mu_Q)\le c(\pi_P)+c(\pi_Q).
\]
Using (5) and (8), we get
\[
f_P+f_Q=c(\mu_P)+c(\mu_Q)\le c(\pi_P)+c(\pi_Q)<e_P+e_Q.
\]
Therefore at least one of the rows \(P,Q\) still has an original \(1\) that has not been used on the small side.

Choose such a row and replace the cut-produced marked \(1\) on its small side by this unused original \(1\). The multiset of the small-side partition of that row is unchanged. Cancel the cut in that row and perform a new cut outside \(P,Q\). The small partition types remain unchanged. After this adjustment, at least one of the two marked small partitions, say \(\mu_I\), is \([1^b]\). If the other one \(\mu_J\) is also \([1^b]\), then
\[
\nu(\mu_I)+\nu(\mu_J)=0\le b-2.
\]
Otherwise, \(\mu_J\) contains the marked part \(1\) together with at least one further part; hence \(|\mu_J|\ge 2\) and therefore
\[
\nu(\mu_J)=b-|\mu_J|\le b-2.
\]
Consequently,
\[
\nu(\mu_I)+\nu(\mu_J)\le 0+(b-2)=b-2.
\]
Thus assertion (5) holds.

For the large side, uncut rows have positive large-side defect. If a cut
row lies outside \(P,Q\), its large-side defect is at least \(\ell-1\),
so the sum of the two cut-row large-side defects is positive. If the two
cuts remain at \(P,Q\), then
\[
\nu(\lambda_P)+\nu(\lambda_Q)=\nu_P+\nu_Q-b-1.
\]
If one of \(P,Q\) satisfies \(e_i<b\), this is positive by \eqref{eq:3}
and \(\nu_j\ge\ell\) for the other. Otherwise \(P,Q\) are the two rows of
maximal defect, and
\[
\nu_P+\nu_Q\ge\frac{2(d+\ell-2)}{\ell+1}>b+1,
\]
where the last inequality follows from
\[
2(\ell b+\ell-2)-(\ell+1)(b+1)=(\ell-1)b+\ell-5>0.
\]
Each actual cut lowers the total defect by \(1\), so
\[
\sum_i\nu(\lambda_i)=d+\ell-2-(b-1)-2=D+\ell-3.
\]
Since the \(\ell-1\) uncut rows have defect at least \(1\),
\[
\nu(\lambda_I)+\nu(\lambda_J)\le D-2.
\]
This proves all assertions.
\end{proof}

\section{Marked gluing of spherical coverings}

\begin{theorem}[Marked gluing]\label{thmB}
Let two connected spherical branched coverings have degrees \(D\) and
\(b\), with monodromy permutations
\[
\alpha_1,\ldots,\alpha_k\in S_A,\qquad
\beta_1,\ldots,\beta_k\in S_B,
\]
where \(|A|=D\), \(|B|=b\), \(A\cap B=\emptyset\), and both products are
the identity. Assume that at the first two branch values, the selected
cycles of \(\alpha_1,\alpha_2\) share a letter \(x\) and have lengths
\(v_1,v_2\), while the selected cycles of \(\beta_1,\beta_2\) share a
letter \(y\) and have lengths \(u_1,u_2\). Then there exists a rational
function of degree \(D+b\) whose first two branch partitions merge the
selected parts into \(v_1+u_1\) and \(v_2+u_2\), all other parts being
the multiset union of the original partitions.
\end{theorem}

\begin{proof}
Let \(\gamma_i=\alpha_i\oplus\beta_i\) and \(\tau=(x\,y)\). Define
\[
\sigma_1=\gamma_1\tau,\qquad
\sigma_2=\tau\gamma_2,\qquad
\sigma_i=\gamma_i\quad(i\ge3).
\]
Then
\[
\sigma_1\cdots\sigma_k
=\gamma_1\tau\tau\gamma_2\cdots\gamma_k
=\gamma_1\cdots\gamma_k=1.
\]
Multiplication by a transposition joining two distinct cycles merges
them, so the first two positions perform the required merges. The new
cycles contain both old supports, so transitivity is preserved on
\(A\sqcup B\). The total defect is
\[
(2D-2)+(2b-2)+2=2(D+b)-2.
\]
Riemann--Hurwitz and the Riemann existence theorem now give the required
rational function.
\end{proof}

\begin{remark}
The common-letter hypothesis is essential: it allows the same
cross-side transposition to be used at both gluing positions. It cannot be
inferred from the separate realizability of the two sides.
\end{remark}

\section{Marked merging lemma and the induction}
The content of the following lemma is taken from \cite{EKS84}; the
marked-part refinement is added for the present application.

\begin{lemma}[Additive defect merging]\label{lem:2}
Let \(\eta,\theta\vdash N\) with
\[
\nu(\eta)+\nu(\theta)\le N-1.
\]
Then there exist permutations \(a,b\) of cycle types \(\eta,\theta\) such
that
\[
\nu(ab)=\nu(\eta)+\nu(\theta).
\]
If one part in each partition is marked, the corresponding cycles can be
made to share a letter, provided the defect sum is at most \(N-2\), or
equals \(N-1\) but the two marked parts are not both \(1\).
When \(N=1\), take the identity permutation.
\end{lemma}

\begin{proof}
We first construct a bipartite forest \(F\) with vertex sets \(U,V\),
where the degrees on \(U\) are given by \(\eta\), the degrees on \(V\) are
given by \(\theta\), and the edge set \(E\) has cardinality \(N\).
Such a forest exists by induction on \(N\). For \(N=1\) it is a single
edge. For \(N\ge2\), put
\[
c=|\eta|+|\theta|-N=N-\nu(\eta)-\nu(\theta)\ge1.
\]
If every vertex on both sides had degree at least \(2\), then
\[
N=\sum_{u\in U}\deg(u)\ge2|\eta|,
\qquad
N=\sum_{v\in V}\deg(v)\ge2|\theta|,
\]
hence \(|\eta|+|\theta|\le N\), contradicting \(c\ge1\). Thus at least
one side has a vertex of degree \(1\). Suppose \(u\in U\) has degree
\(1\). If some \(v\in V\) has degree at least \(2\), delete \(u\) and the
edge \(uv\), and decrease the degree of \(v\) by \(1\). The new degree
sequences \(\eta',\theta'\) have total \(N-1\), and
\[
\nu(\eta')+\nu(\theta')\le (N-1)-1.
\]
By induction there is a forest for \(\eta',\theta'\); adding back the edge
\(uv\) gives the required forest. If every vertex of \(V\) has degree
\(1\), then \(\theta=(1^N)\). In this case each \(u\in U\) is joined to
\(\deg(u)\) distinct vertices of \(V\); since all vertices of \(V\) have
degree \(1\), these stars are disjoint, and we obtain the required
forest directly. The case where a vertex of \(V\) has degree \(1\) is
symmetric. This proves the existence of \(F\).

At each vertex of \(F\), order the incident edges cyclically. Let \(a\)
be the product of the cycles at the vertices of \(U\), and \(b\) the
product of the cycles at the vertices of \(V\), where each edge is
regarded as a letter. Deleting leaf edges successively shows that each
tree of \(F\) contributes exactly one cycle to \(ab\). Hence \(ab\) has
exactly
\[
c=|\eta|+|\theta|-N=N-\nu(\eta)-\nu(\theta)
\]
cycles, so
\[
\nu(ab)=N-c=\nu(\eta)+\nu(\theta).
\]

It remains to treat the marked parts. Let \(u\in U\) and \(v\in V\) be
the marked vertices. We modify \(F\) so that \(uv\) becomes an edge while
preserving the degree sequences and the forest property.

If \(u\) and \(v\) lie in different trees, choose an edge \(uu'\) and an
edge \(vv'\), delete them, and add \(uv\) and \(u'v'\). The two original
trees are split into four subtrees; the new edges join them in pairs, so
the result is again a forest with the same degrees, and now \(uv\) is an
edge.

If \(u\) and \(v\) lie in the same tree and at least one of them, say \(u\),
has degree at least \(2\), choose an edge \(uu'\) not on the \(u\)-\(v\)
path, and choose the edge \(vv'\) on the \(u\)-\(v\) path incident with
\(v\). Delete \(uu'\) and \(vv'\). Since two edges of a tree are removed,
the original tree is split into three subtrees. Add the edges \(uv\) and
\(u'v'\). These two new edges reconnect the three subtrees into a single
tree: the edge \(uv\) joins the component containing \(u\) to the component
containing \(v\), while the edge \(u'v'\) joins the remaining component
containing \(u'\) to the component containing \(v'\), which lies in the
same component as \(u\). No cycle is created, the degree sequences are
preserved, and \(uv\) is an edge of the resulting forest.

If both \(u\) and \(v\) have degree \(1\) and \(c\ge2\), remove the marked
parts \(1\) from \(\eta\) and \(\theta\), i.e.\ set
\[
\eta'=\eta\setminus\{1\},\qquad
\theta'=\theta\setminus\{1\}.
\]
Then \(\eta',\theta'\vdash N-1\), and
\[
\nu(\eta')+\nu(\theta')
=\nu(\eta)+\nu(\theta)
\le N-2=(N-1)-1.
\]
By the existence part applied to \(\eta',\theta'\), the remaining edges
can be completed to a bipartite forest. Together with the isolated edge
\(uv\), this gives the required forest with marked edge \(uv\).

The only excluded case is \(c=1\) with both marked vertices of degree
\(1\); then a connected tree cannot contain the edge \(uv\) without
collapsing to a single edge. This is precisely the case excluded in the
statement.

Finally, suppose the merged data are realized by a transitive permutation
group. Conjugate the whole tuple so that the factor corresponding to
\(ab\) becomes \(ab\) itself; this does not change transitivity. Replace
that factor by the two factors \(a\) and \(b\). The product is unchanged,
and the new generated group contains the old one, so transitivity is
preserved. The marked cycles still share the letter corresponding to the
edge \(uv\). This completes the proof.
\end{proof}

\subsection{Closing the induction for Wei's theorem}

Let \(\pi_1,\ldots,\pi_k\vdash d\) be nontrivial and satisfy
\[
\sum_{i=1}^k\nu(\pi_i)=2d-2,\qquad
k>\ell+1,\qquad
\ell=\min_i|\pi_i|=|\pi_k|.
\]
Assume the base case \(\ell\le2\) is known.

Suppose a counterexample exists, and choose one with \(\ell\) minimal and
then \(k\) minimal. Then \(\ell\ge3\) and
\[
\sum_{i<k}\nu(\pi_i)=d+\ell-2.
\]
If \(k\ge\ell+3\), the two smallest defects among \(\pi_1,\ldots,\pi_{k-1}\)
sum to at most
\[
\frac{2(d+\ell-2)}{k-1}\le d-1.
\]
Lemma~\ref{lem:2} merges them, producing data with shortest length at
most \(\ell\), \(k-1\ge\ell+2\) factors, and the same total defect. By
minimality these data are realizable, and splitting back realizes the
original data, a contradiction. Hence \(k=\ell+2\).

If some \(i<k\) has \(\nu(\pi_i)\le\ell-1\), merging it with \(\pi_k\)
gives a partition of length \(\ell-\nu(\pi_i)\le\ell-1\). The new data
have \(\ell+1\) nontrivial factors and still satisfy the strict
inequality, so minimality gives a contradiction. Thus
\[
\nu(\pi_i)\ge\ell\qquad(i<k).
\]

Choose a smallest part \(b\) of \(\pi_k\), and put \(D=d-b\). Then
\(\ell b\le d\). Apply Theorem~\ref{thmA} to \(\pi_1,\ldots,\pi_{k-1}\),
obtaining \(\mu_i,\lambda_i\) and indices \(I,J\). Set
\[
\mu_k=[b],\qquad \lambda_k=\pi_k\setminus\{b\}.
\]
Then \(\lambda_k\) is nontrivial. On the small side,
\[
\sum_{i<k}\nu(\mu_i)=b-1,\qquad \nu(\mu_k)=b-1,
\]
so the total small-side defect is \(2b-2\). Merge the two marked factors
first using Lemma~\ref{lem:2}; Theorem~\ref{thmA}(5) excludes the marked
exception. Merge the remaining factors successively and take the inverse
of the final \(b\)-cycle as the last factor. This gives a transitive
permutation tuple with total defect \(2b-2\), hence a spherical covering
with the required markings.

On the large side,
\[
\sum_{i<k}\nu(\lambda_i)=D+\ell-3,\qquad
\nu(\lambda_k)=D-\ell+1,
\]
so the total defect is \(2D-2\). Theorem~\ref{thmA} gives
\[
0<\nu(\lambda_I)+\nu(\lambda_J)\le D-2.
\]
Merge these two factors using Lemma~\ref{lem:2}, keeping the marked
cycles intersecting. The auxiliary data have \(\ell+1\) nontrivial
factors, including \(\lambda_k\) of length \(\ell-1\), so their shortest
length is at most \(\ell-1\). Since
\[
\ell+1>(\ell-1)+1,
\]
minimality realizes the auxiliary data. Splitting back gives the large
side with the required markings.

Finally, order the branch values on both sides so that \(I,J\) are
adjacent, construct the marked factor pairs first, and apply
Theorem~\ref{thmB}. The two cut parts are restored as
\(u_I+v_I\) and \(u_J+v_J\), all other parts are restored by multiset
union, and the last row is restored as
\([b]\sqcup\lambda_k=\pi_k\). Thus the original data are realized by a
rational function, contradicting the counterexample assumption.

\section*{Declarations}
\noindent\textbf{Generative AI disclosure.} The authors used OpenAI's ChatGPT as an assisting tool in developing this paper. The authors takes full responsibility for its contents.

\smallskip
\noindent\textbf{Data availability.} This manuscript has no associated data.

\smallskip
\noindent\textbf{Competing interests.} The authors have not disclosed any competing interests.

\section*{Author information}

\noindent
Zhiqiang Wei, Mengyan Xing, Xingya Zhao\\
School of Mathematics and Statistics, Henan University, Zhengzhou, China\\[1mm]
E-mail:
\texttt{10100123@vip.henu.edu.cn} (Zhiqiang Wei);\\
\texttt{104753260851@henu.edu.cn} (Mengyan Xing);\\
\texttt{3383622358zxy@henu.edu.cn} (Xingya Zhao).

\end{document}